# Symplectifying Biarcs

*Geometry of Biarcs and using them with arcs splines*

Stefan Gössner[1]

[1]Dortmund University of Applied Sciences. Department of Mechanical Engineering

## *Abstract*

In this follow-up article to "Symplectification of Circular Arcs and Arc Splines", biarc geometry is examined from a purely geometric point of view. Two given points together with their associated tangent vectors in the plane are sufficient to define two directed, consecutive circular arcs. However, there remains one degree of freedom to determine the join point of both arcs. There are various approaches to this in the literature. A novel one is presented here.

# Content

# 1. Introduction

Piecewise linear or circular curves are commonly used to describe the tool path in CNC Machining and Robot path planning. A planar curve consisting of a number of tangentially joined circular arc and line segments is called an *arc spline* having $G^1$ continuity – which is the condition that the first derivative is continuous [8,10-13]. The second derivative is piecewise constant.

For a pair of oriented $G^1$ Hermite data, two arc segments are generally required to hit both end points and their associated tangents. This pair is called *biarc* [6]. Approximation of smooth planar curves like cubic splines by biarc interpolation has been studied extensively in the past [6,9,11-13]. A purely geometric consideration of a single biarc is discussed in a few papers [6,8,10,12].

The purpose of this publication is to define Biarc geometry using the approach of symplectic geometry in $\mathbb{R}^2$, continuing on from the article on circular arcs [5].

The paper is structured as follows. First, the biarc fundamentals from two sets of Hermite data as boundary condition over the joints circle to different biarc cases are presented. Then, biarc geometry and its one parameter family, are discussed. A new approach for the choice of the joint point location is suggested. Finally, some illustrative examples show, how smoothing of polygonial curves works using different biarc strategies.

## 1.1 Symplectic Geometry in a Nutshell

*Symplectic geometry* in its simplest possible case is the geometry of the plane $\mathbb{R}^2$ [3]. Starting with the *Euclidean vector space* we get the *standard scalar product*, which associates a number to every pair of vectors $\begin{pmatrix} a_1 \\ a_2 \end{pmatrix}$ and $\begin{pmatrix} b_1 \\ b_2 \end{pmatrix}$ in $\mathbb{R}^2$

$$\mathbf{ab} = a_1 b_1 + a_2 b_2 \,.$$

With that the length of a vector and the angle between two vectors is defined. Then we are adding a *complex structure* $\mathbf{J} = \begin{pmatrix} 0 & -1 \\ 1 & 0 \end{pmatrix}$ which, as an *orthogonal operator*, transforms any vector into a *skew-orthogonal* one [1,2]. Now we entered the *complex vector space*.

$$\tilde{\mathbf{a}} = \mathbf{J} \begin{pmatrix} a_1 \\ a_2 \end{pmatrix} = \begin{pmatrix} -a_2 \\ a_1 \end{pmatrix} \,.$$

As a shortcut we will place a *tilde* '~' symbol over the skew-orthogonal vector variable. Yet applying the orthogonal operator to the first vector in the scalar product above gets us to the *skew-scalar product*

$$\tilde{\mathbf{a}}\mathbf{b} = a_1 b_2 - a_2 b_1 \,.$$

The skew-scalar product – called symplectic structure – gives us the *area* of the parallelogram spanned by two vectors $\mathbf{a}$ and $\mathbf{b}$, which is a *directed* or *oriented* area due to its inherent antisymmetry $\tilde{\mathbf{a}}\mathbf{b} = -\mathbf{a}\tilde{\mathbf{b}}$ [1,2,4]. Finally we arrived in the *symplectic vector space*.

The Euclidean, complex and symplectic structure together are named a *compatiple triple*. Having given two of them automatically defines the third. Now we have three compatiple vector spaces – the Euclidean, complex and symplectic vector space in $\mathbb{R}^2$ [1-4].

## 2. Biarc Fundamentals

Two given points $A$ and $B$ together with their associated unit tangent vectors $\mathbf{t}_A$ and $\mathbf{t}_B$ in $\mathbb{R}^2$ are sufficient to define two directed, consecutive arcs. They are starting in $A$ and ending in $B$ while meeting given end point tangents $\mathbf{t}_A$ and $\mathbf{t}_B$.

Both arcs will join at any other given point $\mathbf{J}$ somewhere in the plane. Yet it is a well known fact, that the one-parameter family of points $\mathbf{J}$, in which both arcs meet tangentially, lie on a circle – the *joint circle* (Fig. 1).

Proofs of this can be found in [8,11]. Another one is given below.

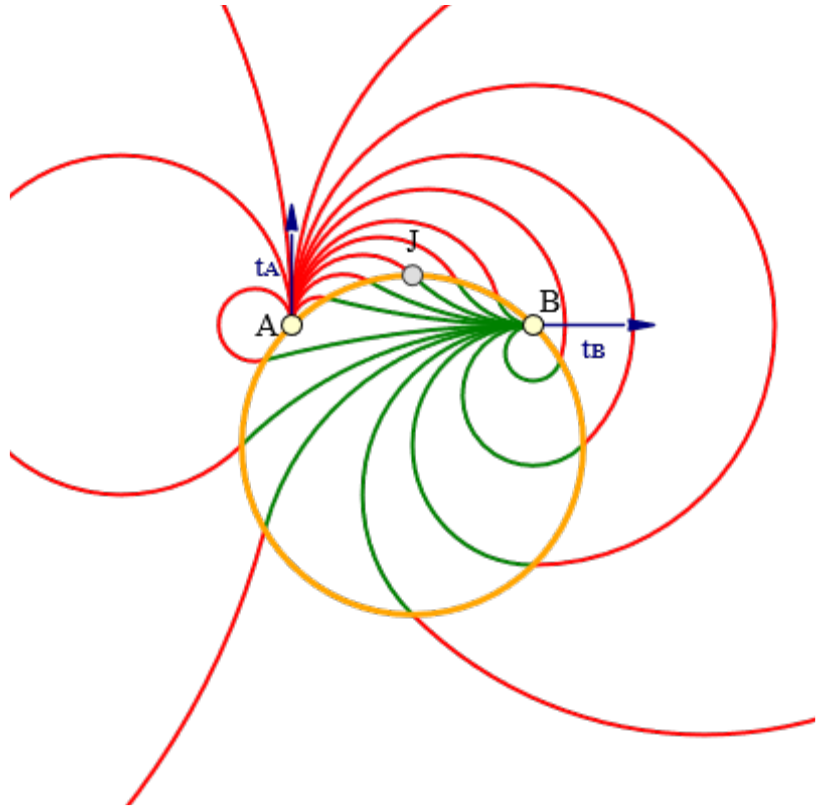

Fig. 1: Some biarcs together with their joint circle.

## 2.1 Two Adjacent Arcs

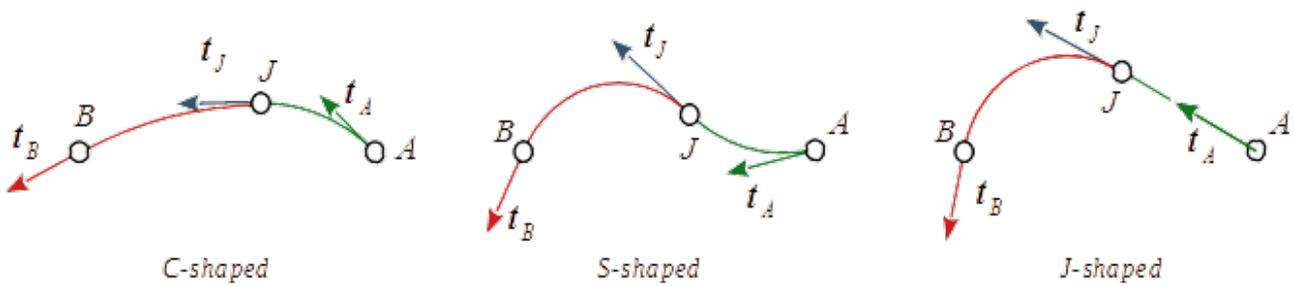

Fig. 2: Different biarc forms.

A *biarc* consists of two circular arcs which satisfy two sets of $G^1$ Hermite data $(\mathbf{A}, \mathbf{t}_A)$ and $(\mathbf{B}, \mathbf{t}_B)$. These two arcs are either circulating in the same direction (*C-shaped*), counter-circulating (*S-shaped*) or one of them has curvature zero and degenerates to a line (*J-shaped*) (Fig. 2) [10,12].

> **Definition 1**
> Two points with associated tangent unit vectors $(\mathbf{A}, \mathbf{t}_A)$ and $(\mathbf{B}, \mathbf{t}_B)$ are interpolated by a *biarc* if and only if
>
> 1. the arc starting in $A$ does so tangential to $t_A$,
> 2. the arc ending in $B$ does so tangential to $t_B$ and
> 3. both arcs have a common tangent $t_J$ in their join point $J$.
>
> These two arcs are called *biarc* [6].

## 2.2 Biarc Angle $\psi$

The angle from tangent vector $\mathbf{t}_A$ to vector $\mathbf{t}_B$ is the characteristic *biarc angle* $\psi$.

$$\sin \psi = \tilde{\mathbf{t}}_A \mathbf{t}_B, \quad \cos \psi = \mathbf{t}_A \mathbf{t}_B, \quad \tan \psi = \frac{\tilde{\mathbf{t}}_A \mathbf{t}_B}{\mathbf{t}_A \mathbf{t}_B} \tag{1}$$

Subsequently we will also use half-angle terms quite frequently.

$$\sin \tfrac{\psi}{2} = \sqrt{\frac{1 - \mathbf{t}_A \mathbf{t}_B}{2}}, \quad \cos \tfrac{\psi}{2} = \sqrt{\frac{1 + \mathbf{t}_A \mathbf{t}_B}{2}}, \quad \tan \tfrac{\psi}{2} = \frac{\tilde{\mathbf{t}}_A \mathbf{t}_B}{1 + \mathbf{t}_A \mathbf{t}_B} = \frac{1 - \mathbf{t}_A \mathbf{t}_B}{\tilde{\mathbf{t}}_A \mathbf{t}_B} \tag{1.1}$$

Please note that angle $\psi$ may deviate by $\pm 2\pi$ in the contexts discussed below.

## 2.3 Joint Circle

Two arcs interpolating two points with associated tangent vectors $(\mathbf{A}, \mathbf{t}_A)$ and $(\mathbf{B}, \mathbf{t}_B)$ might meet in an arbitrary common point $\mathbf{J}$ in the plane, while satisfying conditions 1 and 2 of Definition 1. Adding condition 3 of a common tangent vector $\mathbf{t}_J$ in joint $\mathbf{J}$, results in a one-parameter family of possible join points $\mathbf{J}$ (see Fig. 1).

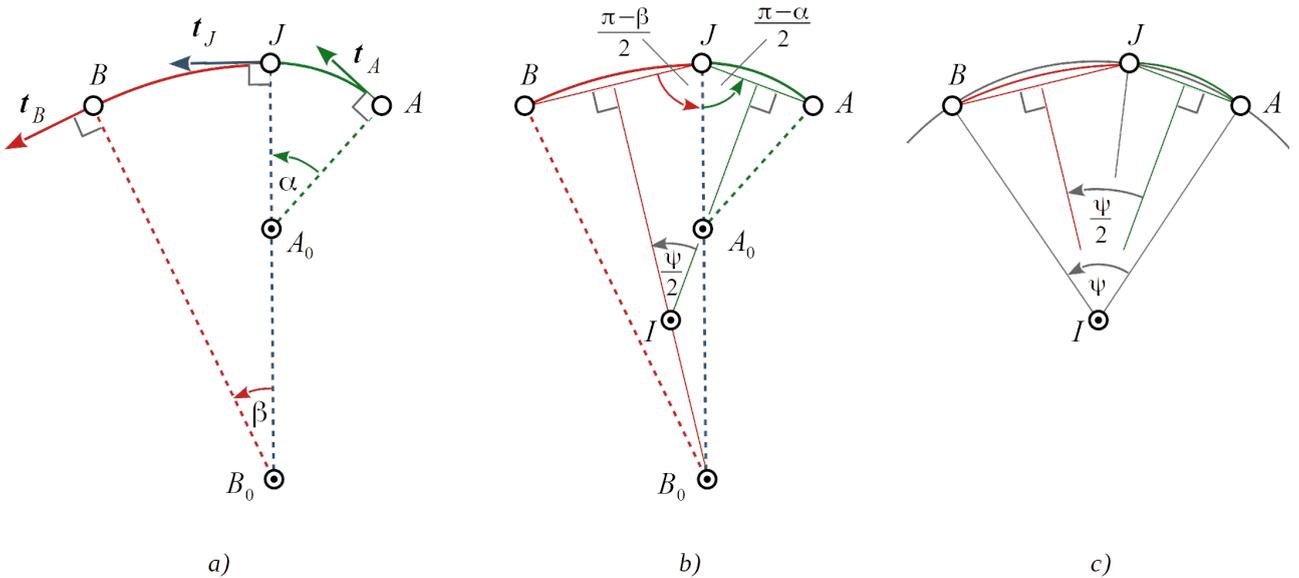

Fig. 3: Geometry of the joint circle.

**Theorem 1:**
The locus of all join points $J$ of a biarc is the *joint circle* through the points $A$ and $B$ with center $I$. Its radial vectors $\mathbf{r}_{IA}$ and $\mathbf{r}_{IB}$ enclose the same directed angle $\psi$ as the tangent vectors $\mathbf{t}_A$ and $\mathbf{t}_B$.

**Proof.**

- Let $\alpha$ and $\beta$ be the directed angles of the arcs starting in $A$ and ending in $B$ (Fig. 3a).
- Denote the intersection point of the perpendicular bisectors of chords $\overline{AJ}$ and $\overline{JB}$ by $I$, which is the center of the circle through the three points $A$, $J$ and $B$ then (Fig. 3b).

- Rotating start vector $\mathbf{t}_A$ about $A_0$ into joint vector $\mathbf{t}_J$ by $\alpha$ and then rotating $\mathbf{t}_J$ about $B_0$ into end vector $\mathbf{t}_B$ by $\beta$ is equivalent to directly rotating vector $\mathbf{t}_A$ about $I$ into vector $\mathbf{t}_B$ by $\psi$, i.e.

$$\psi = \alpha + \beta \,. \tag{2}$$

- Chords $\overline{AJ}$ and $\overline{JB}$ enclose in common point $J$ an angle $\frac{\pi-\alpha}{2} + \frac{\pi-\beta}{2} = \pi - \frac{\alpha+\beta}{2}$ (Fig. 3b). The angle in point $I$ of the quadrilateral opposite to point $J$ is $\frac{\alpha+\beta}{2} = \frac{\psi}{2}$ then, which is a fixed value, thus independent of the location of $J$.
- Due to the symmetry of bisecting isosceles triangles the joint circle center angle $\angle AIB$ equals $\psi$ (Fig. 3c) □

The geometric proof given here is similar to that in [8]. Another proof is to be found in [11].

## 2.4 Joint Circle Radius

> **Theorem 2**:
> The radius of the joint circle through the points $A$ and $B$ with center $I$ is
>
> $$R = \frac{c}{2 \sin \frac{\psi}{2}} \quad for \quad c > 0\,. \tag{3}$$

*Proof.* See Lemma 1 in [5].

Please note, that the radius $R$ is a signed quantity due to equation (3).

## 2.5 Joint Circle Center

> **Theorem 3**:
> The joint circle's center point vector seen from endpoint $A$ is
>
> $$\mathbf{r}_{AI} = \frac{\sin \frac{\psi}{2}\, \mathbf{c} + \cos \frac{\psi}{2}\, \tilde{\mathbf{c}}}{2 \sin \frac{\psi}{2}}\,. \tag{4}$$

*Proof.* See Lemma 2 in [5].

An alternative or trigonometry-free version of equation (4) – using expressions (1.1) – reads

$$\mathbf{r}_{AI} = \frac{\mathbf{c}}{2} + \frac{\tilde{\mathbf{c}}}{2 \tan \frac{\psi}{2}} = \frac{1}{2}\left(\mathbf{c} + \frac{\tilde{\mathbf{t}}_A \mathbf{t}_B}{1 - \mathbf{t}_A \mathbf{t}_B}\tilde{\mathbf{c}}\right)\,. \tag{4.1}$$

## 2.6 Joint Tangent

In the join point $J$ both arcs of a biarc have a common tangent $\mathbf{t}_J$ due to Definition 1.

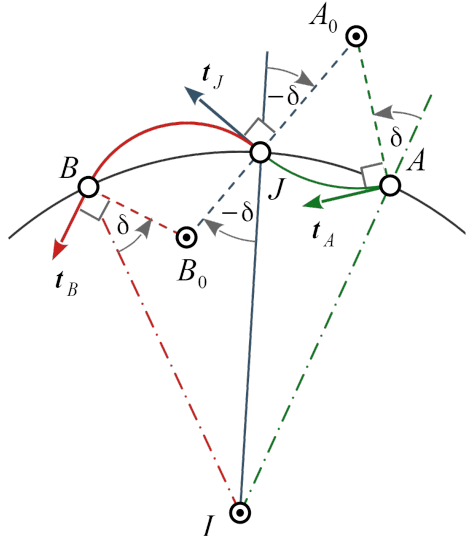

Fig. 4: Join point tangent.

**Theorem 4**: Both end vectors $\mathbf{t}_A$ and $\mathbf{t}_B$ enclose the same directed angle $\delta$ with the joint circle tangent in their corresponding points $A$ and $B$. The tangent vector in the join point $\mathbf{t}_J$ always encloses the negative angle $-\delta$ with the joint circle tangent in $J$.

*Proof.* Due to symmetry of the isosceles triangle $\Delta AA_0J$ the angle from the circle normal in point $J$ is opposite equal to the angle from the circle normal in point $A$ (Fig.4). Same opposite angular relations apply to the isosceles triangle $\Delta JB_0B$. □

## 2.7 Different Biarc Cases

Having given $\mathbf{t}_A$, $\mathbf{t}_B$ and $\mathbf{c}$, we can distinguish between various geometric special cases (Table 1).

In the special case where start and end tangent unit vectors are parallel, i.e. $\mathbf{t}_A = \mathbf{t}_B$, biarc angle $\psi = 0$ due to equation (1) and the joint circle degenerates to a straight line through $A$ and $B$ (cases 5 and 6). If additionally the endpoint vector $\mathbf{c}$ is collinear to both tangent vectors ($\tilde{\mathbf{c}}\mathbf{t}_A = 0$), the biarc also degenerates to a single line (case 6).

If both parallel tangent vectors are not collinear to vector $\mathbf{c}$, i.e. $\mathbf{c}\tilde{\mathbf{t}}_A \neq 0$, valid biarcs are generated by choosing the join point $J$ somewhere on the straight line $\overline{AB}$ (case 5).

If tangent unit vectors are antiparallel ($\psi = \pi$) vector $\mathbf{c}$ is diameter of the joint circle (case 7), even if unit tangents are (anti)collinear to $\mathbf{c}$.

In the general case the end tangent vectors are not parallel, i.e. $\mathbf{t}_A \neq \mathbf{t}_B$. If then $\mathbf{ct}_A > \mathbf{ct}_B$, the end tangent vectors are pointing to the joint circle's inside and the join point tangent points outside (cases 1,4). If $\mathbf{ct}_A < \mathbf{ct}_B$, the end tangent vectors are pointing to the joint circle's outside and the joint tangent points inside (case 3).

In the case of $\mathbf{ct}_A = \mathbf{ct}_B$ all three tangents are tangential to the joint circle and the biarc is reduced to a single arc (case 2).

Table 1: Various geometric situations.

| Case | Geometry | $\mathbf{t}_A = \mathbf{t}_B$ | $\mathbf{ct}_A \gtreqless \mathbf{ct}_B$ | $\tilde{\mathbf{c}}\mathbf{t}_A = 0$ | $\mathbf{t}_A, \mathbf{t}_B \nearrow$ | comment |
|---|---|---|---|---|---|---|
| 1 |  | $\mathbf{t}_A \neq \mathbf{t}_B$ | $\mathbf{ct}_A > \mathbf{ct}_B$ | $\tilde{\mathbf{c}}\mathbf{t}_A \neq 0$ | inside | arc A inner |
| 2 |  | $\mathbf{t}_A \neq \mathbf{t}_B$ | $\mathbf{ct}_A = \mathbf{ct}_B$ | $\tilde{\mathbf{c}}\mathbf{t}_A \neq 0$ | tangential | single arc |
| 3 |  | $\mathbf{t}_A \neq \mathbf{t}_B$ | $\mathbf{ct}_A < \mathbf{ct}_B$ | $\tilde{\mathbf{c}}\mathbf{t}_A \neq 0$ | outside | arc A outer |
| 4 |  | $\mathbf{t}_A \neq \mathbf{t}_B$ | $\mathbf{ct}_A > \mathbf{ct}_B$ | $\tilde{\mathbf{c}}\mathbf{t}_A = 0$ | inside | $\mathbf{t}_A$ collinear |
| 5 |  | $\mathbf{t}_A = \mathbf{t}_B$ | $\mathbf{ct}_A = \mathbf{ct}_B$ | $\tilde{\mathbf{c}}\mathbf{t}_A \neq 0$ | to left | arc A left |
| 6 |  | $\mathbf{t}_A = \mathbf{t}_B$ | $\mathbf{ct}_A = \mathbf{ct}_B$ | $\tilde{\mathbf{c}}\mathbf{t}_A = 0$ | collinear | single line |
| 7 |  | $\mathbf{t}_A = -\mathbf{t}_B$ | $\mathbf{ct}_A = -\mathbf{ct}_B$ | $\tilde{\mathbf{c}}\mathbf{t}_A \neq 0$ | antiparallel | $\mathbf{c}$ is $\emptyset$ |

## 3. Biarc Geometry

Now we have two points $A$ and $B$ with their corresponding tangent unit vectors $t_A$ and $t_B$ as well as the joint circle, on which we are able to choose a join point $J$, and with that we get a particular biarc from its one-parameter family. In section 5 below, we will discuss advantageous ways of selecting the join point $J$. For the moment, we assume that point $J$ has already been chosen. Chord vectors $\mathbf{a} = \mathbf{r}_{AJ}$ and $\mathbf{b} = \mathbf{r}_{JB}$ will be used to describe the biarc geometry (Fig. 5a). The following relation holds

$$\mathbf{a} + \mathbf{b} = \mathbf{c}. \qquad (5)$$

### 3.1 Arc Angles

In order to determine the arc angles we do not use angular arithmetic, nor do we take the directed angle $\alpha$ from the start tangent to the joint tangent vector and $\beta$ from the joint tangent to the end tangent vector.

Both methods suffer from a limited angular range $\alpha, \beta \in [-\pi, \pi]$. We rather take a half angle approach.

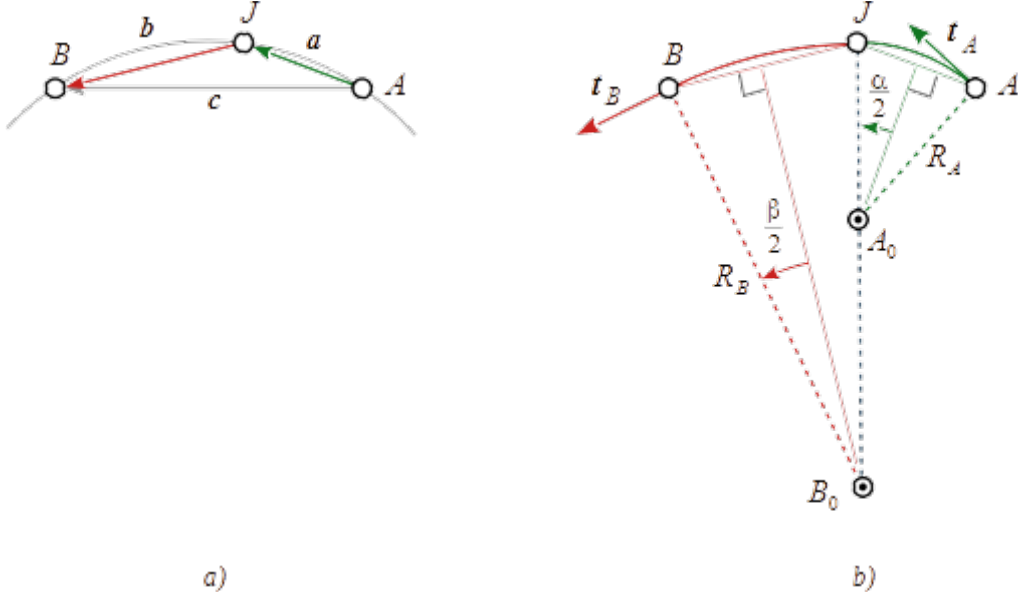

Fig. 5: Chord vectors, arc angles and radii.

**Theorem 5**:
Half arc angle $\frac{\alpha}{2}$ is measured from the start tangent $\mathbf{t}_A$ to chord vector $\mathbf{a}$ and $\frac{\beta}{2}$ from chord vector $\mathbf{b}$ to end tangent vector $\mathbf{t}_B$.

$$\tan \frac{\alpha}{2} = \frac{\tilde{\mathbf{t}}_A \mathbf{a}}{\mathbf{t}_A \mathbf{a}} \quad and \quad \tan \frac{\beta}{2} = \frac{\tilde{\mathbf{b}} \mathbf{t}_B}{\mathbf{b} \mathbf{t}_B} \quad with \quad \alpha, \beta \in [-2\pi, 2\pi] \tag{6}$$

Arc angles $\alpha$ and $\beta$ are positive with arcs running counterclockwise, otherwise negative.

*Proof.* Based on the similarity of right-angled triangles (Fig. 5b). □

Note, that the sum of both $\alpha$ and $\beta$ is constrained by equation (2). Note again, that the sum may differ by $\pm 2\pi$ from $\psi$.

## 3.2 Arc Radii

**Theorem 6**:
Signed Arc radii for a given biarc instance, expressed by tangent vectors $\mathbf{t}_A$ and $\mathbf{t}_B$ and chord vectors $\mathbf{a}$ and $\mathbf{b}$, are

$$R_A = \frac{a^2}{2\tilde{\mathbf{t}}_A \mathbf{a}} = \frac{a}{2\sin\frac{\alpha}{2}} \quad and \quad R_B = \frac{b^2}{2\tilde{\mathbf{b}} \mathbf{t}_B} = \frac{b}{2\sin\frac{\beta}{2}}. \tag{7}$$

Radii values $R_A$ and $R_B$ are positive with arcs running counterclockwise, otherwise negative.

*Proof.* Loop closure equation of half isosceles triangle $\Delta AA_0J$ gives $R_A \tilde{\mathbf{t}}_A - \lambda \tilde{\mathbf{a}} - \frac{1}{2}\mathbf{a} = \mathbf{0}$. Multiplication by vector $\mathbf{a}$ eliminates the second summand and allows to resolve for $R_A$. Applying the same procedure to half isosceles triangle $\Delta JB_0B$ leads to the arc radii (7). □

Arc centers $A_0$ and $B_0$ are to be found from $A$ and $B$ along their normal unit vectors $\tilde{\mathbf{t}}_A$ and $\tilde{\mathbf{t}}_B$ (Fig. 5b).

$$\mathbf{r}_{AA_0} = R_A \tilde{\mathbf{t}}_A \quad and \quad \mathbf{r}_{BB_0} = R_B \tilde{\mathbf{t}}_B \tag{7.1}$$

# 4. One-parameter Family of Biarcs

For two sets of $G^1$ Hermite data $(\mathbf{A}, \mathbf{t}_A)$ and $(\mathbf{B}, \mathbf{t}_B)$ a one-parameter family of biarcs exist. The locus of all possible join points is the joint circle. So a parameterization via a circular angle seems to be beneficial and convenient for the choice of the joint location.

## 4.1 Joint Circle Origin

Let $\varphi$ be the angular parameter for the location of the join point on the joint circle. In $\mathbb{R}^2$ the standard way to measure the (absolute) angle of a single vector is done with respect to the x-axis rightward. This is not, what we want to adopt. We rather prefer a biarc-specific joint location origin. So we take the circle point on the perpendicular bisector of the line $\overline{AB}$ nearest to it as the origin $J_0$ with $\varphi = 0$. Angular parameter $\varphi$ is measured from here anticlockwise (Fig. 6).

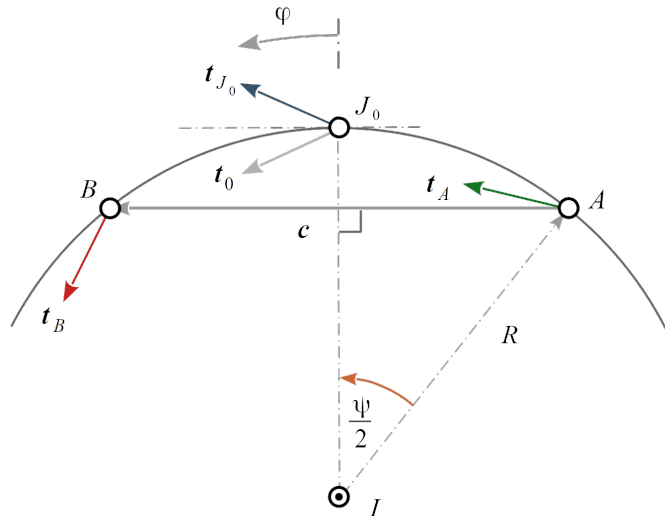

Fig. 6: Joint circle origin location.

**Lemma 1**: The biarc's joint circle origin is located at (seen from joint circle center $I$)

$$\mathbf{r}_{IJ_0} = -\frac{\tilde{\mathbf{c}}}{2 \sin \frac{\psi}{2}} . \tag{8}$$

or seen from start point $A$

$$\mathbf{r}_{AJ_0} = \frac{1}{2}\left(\mathbf{c} - \tan\frac{\psi}{4}\tilde{\mathbf{c}}\right) . \tag{8.1}$$

*Proof.* We get $\mathbf{r}_{IJ_0}$ via rotating $\mathbf{r}_{IA}$ by $\frac{\psi}{2}$, i.e. $\mathbf{r}_{IJ_0} = \cos\frac{\psi}{2}\mathbf{r}_{IA} + \sin\frac{\psi}{2}\tilde{\mathbf{r}}_{IA}$. Inserting $\mathbf{r}_{IA}$ from expression (4) leads to equation (8). Using $\mathbf{r}_{AJ_0} = \mathbf{r}_{AI} + \mathbf{r}_{IJ_0}$ with equations (4) and (8) together with expression $\tan\frac{x}{4} = \frac{1-\cos\frac{x}{2}}{\sin\frac{x}{2}}$ yields equation (8.1). □

## 4.2 Joint Circle Origin Tangent

**Lemma 2**: The biarc's joint circle origin tangent is

$$\mathbf{t}_{J_0} = \cos\frac{\psi}{2}\mathbf{t}_A^* - \sin\frac{\psi}{2}\tilde{\mathbf{t}}_A^* \quad with \quad \mathbf{t}_A^* = \frac{(\mathbf{ct}_A)\mathbf{c} - (\tilde{\mathbf{c}}\mathbf{t}_A)\tilde{\mathbf{c}}}{c^2} \tag{9}$$

*Proof.* We can get the origin tangent $\mathbf{t}_{J_0}$ by one of two methods:

1. Rotate tangent $\mathbf{t}_A$ by $\frac{\psi}{2}$, then reflect it at $\mathbf{c}$.
2. Reflect tangent $\mathbf{t}_A$ at $\mathbf{c}$, then rotate it by $-\frac{\psi}{2}$.

We follow the second method and reflect tangent $\mathbf{t}_A$ at $\mathbf{c}$ first and get $\mathbf{t}_A^*$. This conforms to Grassmann identity (I.2) with its sign of the second summand negated [4]. Rotate $\mathbf{t}_A^*$ then by $-\frac{\psi}{2}$ results in equation (9). □

## 4.3 Parametric Joint Location

An angle $\varphi$ as a position parameter of the join point has a certain disadvantage in the practically significant special case $\psi = 0$. We eliminate this shortcoming by parametrize $\varphi$ itself as $\varphi = u\frac{\psi}{2}$. For the same reason, we describe the position of $J$ not from the center of the circle $I$, but rather from the biarc start point $A$.

> **Theorem 7**:
> The chord vector $\mathbf{a}$ from the biarc start point $A$ to join point $J$ as a function of its position parameter $u$ is
>
> $$\mathbf{a}(u) = \begin{cases} \dfrac{(\sin\frac{\psi}{2} + \sin(u\frac{\psi}{2}))\,\mathbf{c} + (\cos\frac{\psi}{2} - \cos(u\frac{\psi}{2}))\,\tilde{\mathbf{c}}}{2\sin\frac{\psi}{2}} & \text{if } \psi \neq 0 \\ \dfrac{1+u}{2}\mathbf{c} & \text{if } \psi = 0 \end{cases} \qquad (10)$$

*Proof.* In a first step we apply a rotate transformation to the origin vector $\mathbf{r}_{IJ_0}$ by angular parameter $\varphi$ and use expression (8).

$$\mathbf{r}_{IJ}(\varphi) = \cos\varphi\,\mathbf{r}_{IJ_0} + \sin\varphi\,\tilde{\mathbf{r}}_{IJ_0} = \frac{\sin\varphi\,\mathbf{c} - \cos\varphi\,\tilde{\mathbf{c}}}{2\sin\frac{\psi}{2}}$$

For numerical reasons, we want to avoid direct dependencies on the joint circle center point position $I$. Getting the chord vector $\mathbf{a} = \mathbf{r}_{AJ} = \mathbf{r}_{AI} + \mathbf{r}_{IJ}$ by expression (4) and favorable parametrization $\varphi = u\frac{\psi}{2}$ obtains the non-zero case of equation (10)

$$\mathbf{a}(u) = \frac{(\sin\frac{\psi}{2} + \sin(u\frac{\psi}{2}))\,\mathbf{c} + (\cos\frac{\psi}{2} - \cos(u\frac{\psi}{2}))\,\tilde{\mathbf{c}}}{2\sin\frac{\psi}{2}}\,.$$

In case of $\psi = 0$ an indeterminate expression $(\frac{0}{0})$ results. Applying l'Hospital's rule leads to

$$\lim_{\frac{\psi}{2}\to 0}\mathbf{r}_{AJ} = \lim_{\frac{\psi}{2}\to 0}\frac{(\cos\frac{\psi}{2} + u\cos(u\frac{\psi}{2}))\,\mathbf{c} - (\sin\frac{\psi}{2} - u\sin(u\frac{\psi}{2}))\,\tilde{\mathbf{c}}}{2\cos\frac{\psi}{2}} = \frac{1+u}{2}\mathbf{c}\,,$$

a linear function, which conforms to the fact that the joint circle degenerates to a straight line in that case. □

You might confirm quite easily, that $\mathbf{a}(-1) = \mathbf{0}$, $\mathbf{a}(0) = \mathbf{r}_{AJ_0}$ and $\mathbf{a}(1) = \mathbf{r}_{AB}$.

In the practically relevant case, where chord vector $\mathbf{a}$ is given, the associated position parameter $u$ is obtained from

$$u = \begin{cases} 2\frac{\varphi}{\psi} \text{ with } \tan\varphi = \dfrac{2\mathbf{ac} - c^2}{2\tilde{\mathbf{a}}\mathbf{c} - \frac{c^2}{\tan\frac{\psi}{2}}} & \text{if } \psi \neq 0 \\ 2\dfrac{\mathbf{ac}}{c^2} - 1 & \text{if } \psi = 0 \end{cases} \qquad (10.1)$$

## 4.4 Parametric Joint Tangent

> **Theorem 8**:
> The biarc tangent in join point $J$ as a function of parameter $u$ is
> 
> $$\mathbf{t}_J(u) = \begin{cases} \cos(u\frac{\psi}{2})\,\mathbf{t}_{J_0} + \sin(u\frac{\psi}{2})\,\tilde{\mathbf{t}}_{J_0} & \text{if } \psi \neq 0 \\ \mathbf{t}_A^* & \text{if } \psi = 0 \end{cases} \qquad (11)$$

*Proof.* We apply a rotate transformation to the origin's tangent vector $\mathbf{t}_{J_0}$ (12) by $\varphi = u\frac{\psi}{2}$. In case of $\psi = 0$ joint tangent vector $\mathbf{t}_J$ degenerates to $\mathbf{t}_{J_0}$, which itself degenerates to constant vector $\mathbf{t}_A^*$ (see equation (9)).  □

# 5. Choice of the Joint Parameter

Various general approaches can be found in the literature for selecting a specific biarc from its one-parameter family. Here we discuss

- *Equal chord* biarc
- *Parallel tangent* biarc
- *J-shaped* biarc

and propose a new strategy

- *Cubic midpoint* biarc

These approaches mostly lead directly to the joint circle parameter $\varphi$.

## 5.1 Equal Chord Biarc

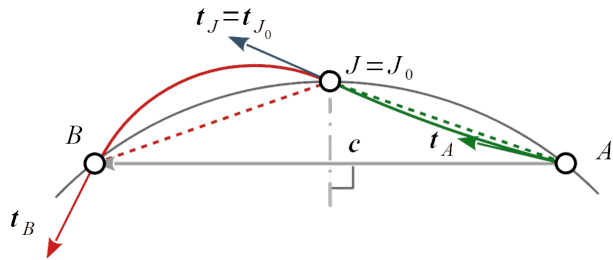

Fig. 7: Equal chord biarc.

The "*equal chord*" biarc is a popular proposal in the literature on biarc approximation [11]. Here the chords $\overline{AJ}$ and $\overline{JB}$ have equal length (Fig. 7). Not surprisingly, this is the trivial case $\varphi = u = 0$ with $\mathbf{J} = \mathbf{J}_0$. This is a numerically robust approach that also works for parallel end tangents. It usually leads to pleasing results for biarc splines.

## 5.2 Parallel Tangent Biarc

The *parallel tangent* biarc, where the joint tangent $\mathbf{t}_J$ is directed parallel to line $\overline{AB}$, is frequently found in biarc literature. It only works for non-parallel end tangents ($\psi \neq 0$).

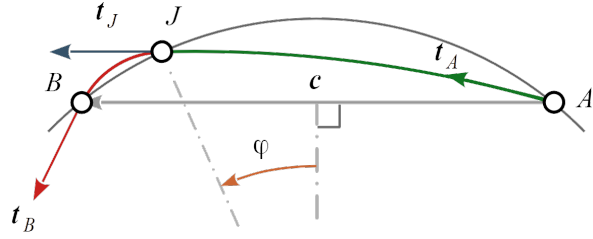

Fig. 8: Parallel tangent biarc.

**Theorem 9**: Parallel tangent biarc

The joint tangent vector is directed parallel to **c** if and only if

$$\tan \varphi = \frac{\tilde{\mathbf{t}}_{J_0}\mathbf{c}}{\mathbf{t}_{J_0}\mathbf{c}} \tag{12}$$

*Proof.* Requesting the joint tangent vector to be parallel to **c**, i.e. $\tilde{\mathbf{c}}\mathbf{t}_J = 0$ (Fig. 8) using the joint tangent (11) from Theorem 8 gives

$$\tilde{\mathbf{c}}(\cos\varphi \mathbf{t}_{J_0} + \sin\varphi \tilde{\mathbf{t}}_{J_0}) = 0$$

Resolving for $\tan\varphi$ leads to expression (12). □

This *parallel tangent* biarc is not as robust as the *equal chord* biarc. Two solutions exist; one for parallel $\mathbf{t}_J \leftarrows \mathbf{c}$ and another one for antiparallel $\mathbf{t}_J \leftrightarrows \mathbf{c}$ vectors. Smooth approximation is expected only if $\varphi \in [-\frac{\psi}{2}, \frac{\psi}{2}]$.

## 5.3 J-shaped Biarc

If start and end tangent point inwards to the joint circle, the starting arc has the potential to be a line (case 1 in Table 1). If both tangents are pointing outside from the joint circle, the ending arc might be a line (case 3 in Table 1).

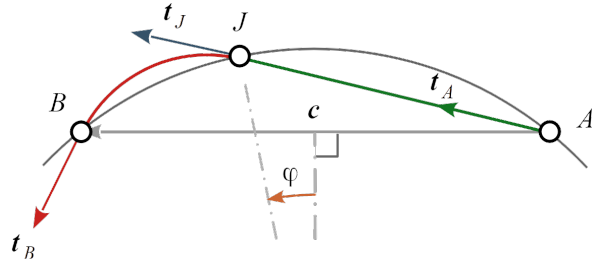

Fig. 9: J-shaped biarc.

**Theorem 10**: J-shaped biarc

Either the starting arc or the ending arc can be forced to be a line.

$$\tan\varphi = \begin{cases} \frac{\tilde{\mathbf{t}}_{J_0}\mathbf{t}_A}{\mathbf{t}_{J_0}\mathbf{t}_A} & \text{if } \mathbf{ct}_A > \mathbf{ct}_B \quad \text{starting arc} \\ 0 & \text{if } \mathbf{ct}_A = \mathbf{ct}_B \quad \text{single arc} \\ \frac{\tilde{\mathbf{t}}_{J_0}\mathbf{t}_B}{\mathbf{t}_{J_0}\mathbf{t}_B} & \text{if } \mathbf{ct}_A < \mathbf{ct}_B \quad \text{ending arc} \end{cases} \tag{13}$$

*Proof.* The starting arc degenerates to a line, if $\mathbf{t}_A$ and $\mathbf{t}_J$ are collinear, i.e. $\tilde{\mathbf{t}}_A\mathbf{t}_J = 0$. Using the joint tangent (11) gives

$$\tilde{\mathbf{t}}_A(\cos\varphi \mathbf{t}_{J_0} + \sin\varphi \tilde{\mathbf{t}}_{J_0}) = 0$$

Resolving for $\tan\varphi$ leads to the corresponding case of expression (13). The same applies to the ending arc, if $\mathbf{t}_J$ and $\mathbf{t}_B$ are collinear, i.e. $\tilde{\mathbf{t}}_J\mathbf{t}_B = 0$. □

Similar to the *parallel tangent* biarc is a smooth approximation expected only if $\varphi \in [-\frac{\psi}{2}, \frac{\psi}{2}]$ with the *J-shaped* biarc.

## 5.4 Curvature Constrained Biarc

Suppose we want to define the curvature or radius of one of the two arcs in order to

- adopt the radius of the neighboring arc or
- not to fall below a minimum radius amount.

This is not a general biarc strategy, but an additional method used under certain circumstances.

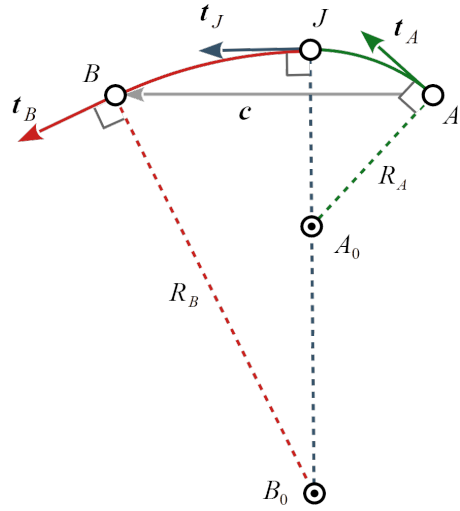

Fig. 10: Curvature constrained biarc.

In this case one radius may be defined and all other missing biarc parameters can be derived. We start with the vector loop closure $AA_0JB_0B$ (Fig. 15), i.e. $R_A\tilde{\mathbf{t}}_A - R_A\tilde{\mathbf{t}}_J + R_B\tilde{\mathbf{t}}_J - R_B\tilde{\mathbf{t}}_B - \mathbf{c} = \mathbf{0}$ or

$$(R_A - R_B)\tilde{\mathbf{t}}_J = R_A\tilde{\mathbf{t}}_A - R_B\tilde{\mathbf{t}}_B - \mathbf{c}\,. \tag{14.1}$$

Squaring eliminates the unknown vector $\mathbf{t}_J$, outmultiplying and combining yields

$$R_A R_B(\mathbf{t}_A\mathbf{t}_B - 1) + R_A\tilde{\mathbf{t}}_A\mathbf{c} - R_B\tilde{\mathbf{t}}_B\mathbf{c} = \frac{c^2}{2}$$

Resolving for the unknown radius of $R_A$ or $R_B$ results in

$$R_A = \frac{\frac{c^2}{2} + R_B\tilde{\mathbf{t}}_B\mathbf{c}}{R_B(\mathbf{t}_A\mathbf{t}_B - 1) + \tilde{\mathbf{t}}_A\mathbf{c}} \quad \text{and} \quad R_B = \frac{\frac{c^2}{2} - R_A\tilde{\mathbf{t}}_A\mathbf{c}}{R_A(\mathbf{t}_A\mathbf{t}_B - 1) - \tilde{\mathbf{t}}_B\mathbf{c}} \tag{14}$$

Now that we know both radii $R_A$ and $R_B$, we can take equation (14.1) and solve for $\mathbf{t}_J$

$$\mathbf{t}_J = \frac{R_A\tilde{\mathbf{t}}_A - R_B\tilde{\mathbf{t}}_B - \mathbf{c}}{R_A - R_B}\,, \tag{14.2}$$

and then by equation (14.2) we get the corresponding join point location via

$$\tan\varphi = \frac{\tilde{\mathbf{t}}_J \mathbf{t}_{J_0}}{\mathbf{t}_J \mathbf{t}_{J_0}}\,. \tag{14.3}$$

Be aware of the fact, that you not only might choose the amount of the radius, but also need to select the correct sign in order to obey the rule of positive arc length.

## 5.5 Cubic Midpoint Biarc

The idea is to define a cubic Bezier curve corresponding to the given Hermite datasets $(\mathbf{A}, \mathbf{t}_A)$ and $(\mathbf{B}, \mathbf{t}_B)$ with two equidistant control points $\mathbf{A}_1$ and $\mathbf{B}_1$, where their distance is chosen so that the center of the Bezier curve $M$ lies on the joint circle. This is a pragmatic approach inspired by papers in which existing curves are to be approximated by biarcs [11].

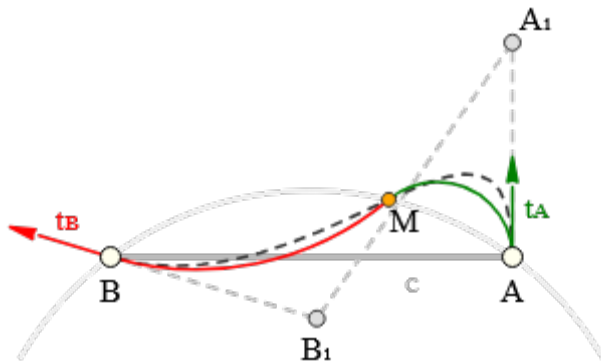

Fig. 11: Midpoint of cubic Bezier chosen as the biarc join point.

> **Theorem 11**:
> A cubic Bezier curve $\mathbf{Q}(t)$ matching two sets of Hermite data $(\mathbf{A}, \mathbf{t}_A)$ and $(\mathbf{B}, \mathbf{t}_B)$ with two equidistant control points $\mathbf{A}_1 = \mathbf{A} + h\mathbf{t}_A$ and $\mathbf{B}_1 = \mathbf{B} - h\mathbf{t}_B$ intersects the joint circle of the respective biarc in a single point $J$. Factor $h$ can be selected so that the midpoint of the cubic $\mathbf{Q}(\frac{1}{2})$ lies on the joint circle in order to become the join point $J$ of the associated biarc.
>
> $$\mathbf{a} = \tfrac{1}{2}\mathbf{c} + \tfrac{3}{8}h\mathbf{t}_{BA} \tag{15}$$
>
> $$\text{with} \quad \mathbf{t}_{BA} = \mathbf{t}_A - \mathbf{t}_B, \quad \kappa = \frac{\tilde{\mathbf{t}}_{BA}\mathbf{c}}{\mathbf{t}_{BA}^2}, \quad h = \frac{-4\kappa}{3\tan\frac{\psi}{2}} + \sqrt{\frac{16\kappa^2}{9\tan^2\frac{\psi}{2}} + \frac{4c^2}{9\sin^2\frac{\psi}{2}}}$$

*Proof.* We define a cubic Bezier curve $\mathbf{Q}(t)$ matching biarc end points $A, B$ and tangents $\mathbf{t}_A, \mathbf{t}_B$ (Fig. 11).

$$\mathbf{Q}(t) = (1-t)^3\,\mathbf{A} + 3(1-t)^2 t\,\mathbf{A}_1 + 3(1-t)t^2\,\mathbf{B}_1 + t^3\,\mathbf{B}\,.$$

Control points $\mathbf{A}_1 = \mathbf{A} + h\mathbf{t}_A$ and $\mathbf{B}_1 = \mathbf{B} - h\mathbf{t}_B$ are located equidistantly to their end points. Taking $\mathbf{B} = \mathbf{A} + \mathbf{c}$ into account yields

$$\mathbf{Q}(t) = (1-t)^3\,\mathbf{A} + 3(1-t)^2 t\,(\mathbf{A} + h\mathbf{t}_A) + 3(1-t)t^2\,(\mathbf{A} + \mathbf{c} - h\mathbf{t}_B) + t^3\,(\mathbf{A} + \mathbf{c})\,.$$

Outmultiplying and simplifying gives

$$\mathbf{Q}(t) = \mathbf{A} + (3-2t)t^2\,\mathbf{c} + 3t(1-t)^2\,h\mathbf{t}_A - 3t^2(1-t)\,h\mathbf{t}_B\,.$$

The curve's midpoint is located at $\mathbf{M} = \mathbf{Q}(\tfrac{1}{2}) = \mathbf{A} + \tfrac{1}{2}\mathbf{c} + \tfrac{3}{8}h(\mathbf{t}_A - \mathbf{t}_B)$, so with $\mathbf{t}_{BA} = \mathbf{t}_A - \mathbf{t}_B$ we get vector

$$\mathbf{r}_{AM} = \tfrac{1}{2}\mathbf{c} + \tfrac{3}{8}h\mathbf{t}_{BA}\,.$$

In order to place the cubic midpoint $M$ onto the joint circle we require $M$ to satisfy Thales' theorem, which in vectorial notation reads: $2\mathbf{r}_{AI} = \mathbf{r}_{AM} + \tau \tilde{\mathbf{r}}_{AM}$. Multiplying that by $\mathbf{r}_{AM}$ eliminates the last term containing unknown scalar $\tau$, i.e. $2\mathbf{r}_{AI}\,\mathbf{r}_{AM} = r_{AM}^2$. Using equation (4.1) for $\mathbf{r}_{AI}$ yields

$$\tfrac{1}{2}\left(\mathbf{c} + \frac{\tilde{\mathbf{c}}}{\tan\tfrac{\psi}{2}}\right)\left(\mathbf{c} + \tfrac{3}{4}h\mathbf{t}_{BA}\right) = \tfrac{1}{4}\left(\mathbf{c} + \tfrac{3}{4}h\mathbf{t}_{BA}\right)^2 .$$

Outmultiplying and simplifying by expression (3.2) gives us a quadratic equation in $h$

$$h^2 + \frac{8\kappa}{3\tan\tfrac{\psi}{2}} h - \frac{4c^2}{9\sin^2\tfrac{\psi}{2}} = 0 \quad with \quad \kappa = \frac{\tilde{\mathbf{t}}_{BA}\mathbf{c}}{\mathbf{t}_{BA}^2}$$

with its positive solution

$$h = \frac{-4\kappa}{3\tan\tfrac{\psi}{2}} + \sqrt{\frac{16\kappa^2}{9\tan^2\tfrac{\psi}{2}} + \frac{4c^2}{9\sin^2\tfrac{\psi}{2}}} .$$

Midpoint $M$ is joint $J$ of the biarc then. □

Incidentally, $\kappa$ has a geometric meaning, the explanation of which is beyond the scope of this article. The joint parameter $u$ corresponding to chord vector $\mathbf{a}$, given by equation (15), is obtained from equation (10.1).

# 6. Examples

## 6.1 Comparing Biarcs

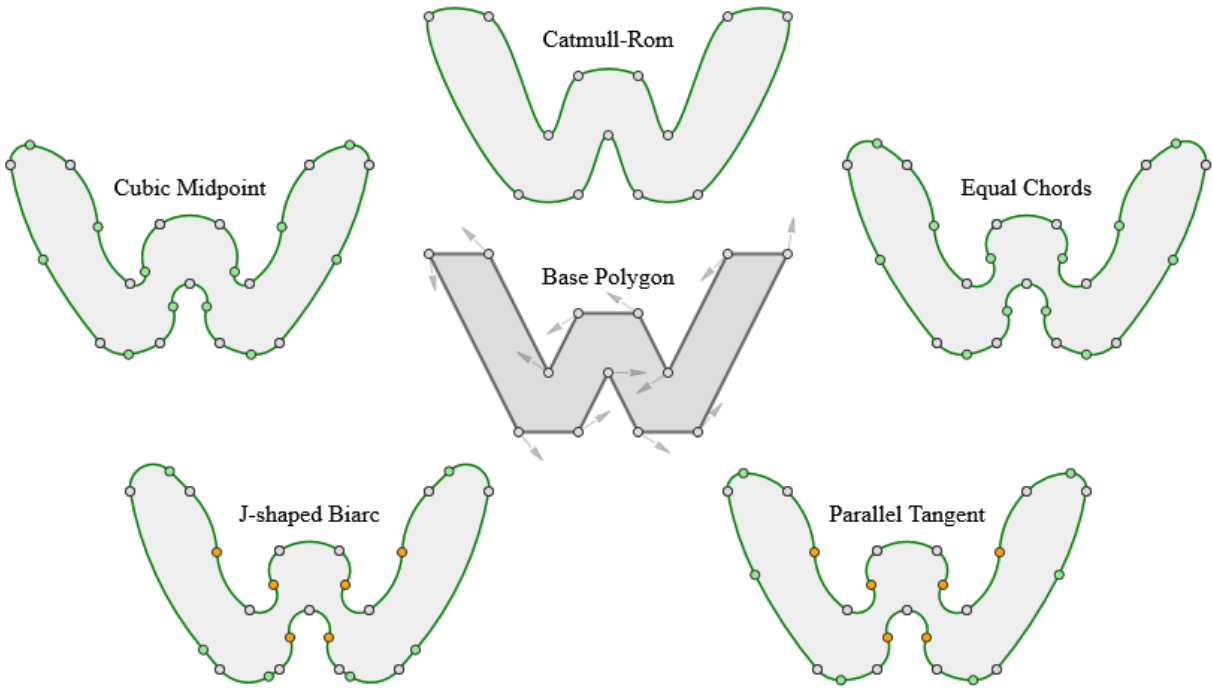

Fig. 12: Arc splines based on biarcs with different joint parameter choice methods.

This first example uses the contour of the letter 'W' as a 13-sided base polygon (central in Fig. 12). Due to its acute angles, it can be considered a stress test for Biarc splines. Simple arc splines, which are discussed in [5], are unable to achieve an appealing result due to that angles.

The curve in Figure 12 on top is a *Centripetal Catmull-Rom Spline* [14], here with the role of a reference curve. The tangents shown in the vertices of the central base polygon are taken from the *Catmull-Rom* algorithm, according which the tangent direction in a vertex is parallel to the line from its predecessor to its successor vertex. The tangents are used to generate the $G^1$ Hermite data sets for the four biarc curves shown.

The two lower curves use the *J-shaped Biarc* and the *Parallel Tangent Biarc*. The light dots are the polygon verticies and the green dots are the joints of the respective biarcs. A closer inspection shows, that one of the green dot's neighboring arcs is a straight line on the left (J-shape) and the connection line of the green dot's neighboring verticies is parallel to the green dot's tangent. Some biarcs of both curves wouldn't have been able to generate a smooth approximation (see section 5.2 and 5.3). So the more robust *Equal Chord Biarc* has been used instead, whose joints were marked orange then.

The two upper curves use the *Cubic Midpoint Biarc* and the *Equal Chords Biarc*. Both handle the critical sharp angles in a robust manner and show visually pleasing results. The newly introduced *Cubic Midpoint Biarc* approximates the sharp angles better than the *Equal Chords Biarc*, thus approximates the reference curve visually best.

## 6.2 Comparing Biarcs and Arc Splines

The second example shows the outline of an animal (Fig. 13). The central 18-sided base polygon is equipped again with tangents at the vertices according to the Catmull-Rom algorithm. In contrast to example 6.1 with the letter "W" no polygon angle herein is significantly sharper than 90°. This has a considerable positiv impact on the quality of the Biarc approximation.

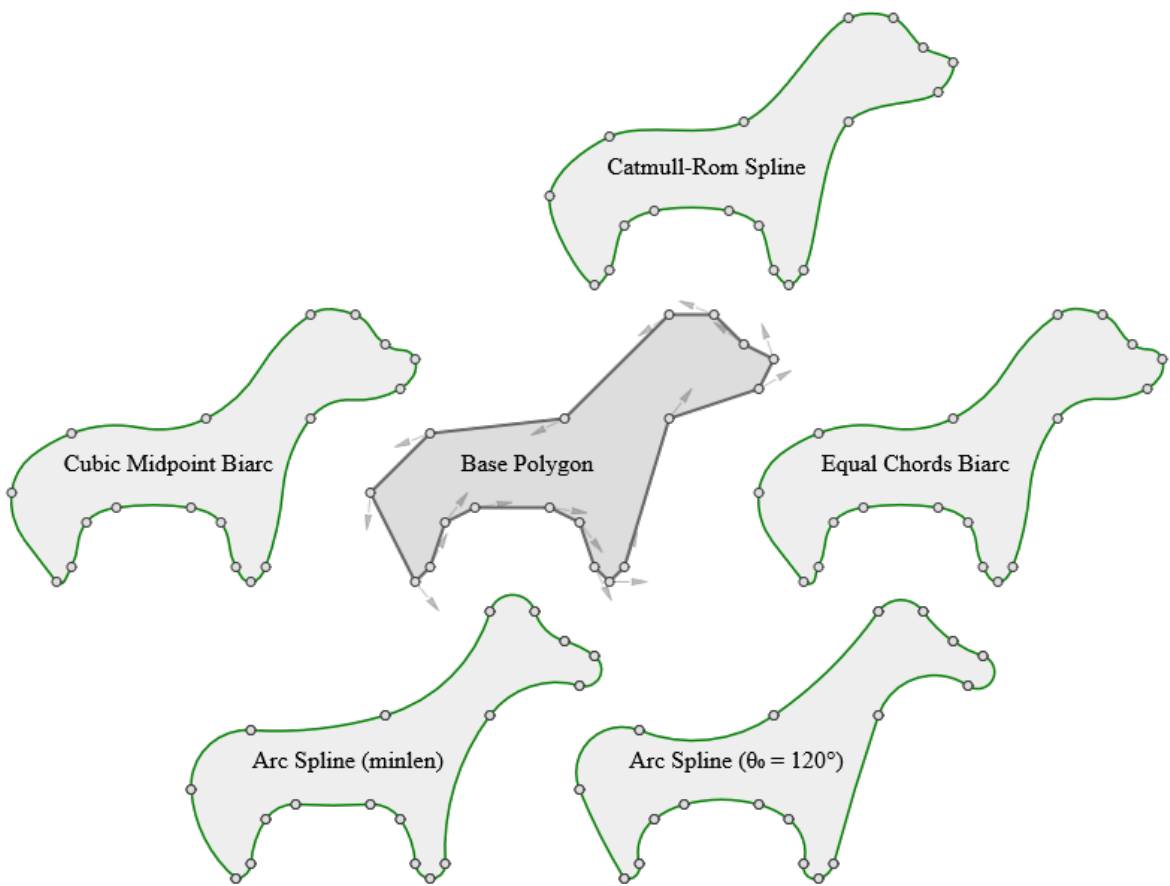

Fig. 13: Biarcs versus simple arcs.

Visually, the Biarc splines on the right and left are barely distinguishable from the upper centripetal Catmull-Rom reference curve. Note also that the *Equal Chords* (right) and the *Cubic Midpoint* Biarc (left) are almost identical. Interestingly, the *Parallel Tangent* and the *J-shaped* Biarc (not shown in Fig. 13) couldn't completely approximate this example on their own, but had to be replaced to a substantial extent ($\approx 50\%$) by the *Equal Chords* Biarc again.

The lower pair of curves are simple arc splines [5]. The right one is based on the freely selectable starting arc angle $\theta_0 = 120°$, the left one is minimized according to the total curve length. The quality of that lower left approximation is not as bad as it seems at first glance, considering that arc splines have as few circular arcs as polygons have sides, while biarc splines have twice as many.

For smooth curve approximation, such as that used in the design of cam and follower mechanisms or in the design of robot paths with the design criterion "as few arcs as possible", the use of simple arc splines may be a serious consideration.

## 6.3 Manual Calculation Example

An example of manual calculating a single biarc shall illustrate the practical use of the equations presented so far. It also shows that the effort required to calculate a biarc with paper and pencil is quite manageable.

Figure 14 shows two given sets of $G^1$ Hermite data $(\mathbf{A}, \mathbf{t}_A)$ and $(\mathbf{B}, \mathbf{t}_B)$. Find the *Equal Chords* and The *Cubic Midpoint* Biarc from this data. An abstract length unit $e$ is used.

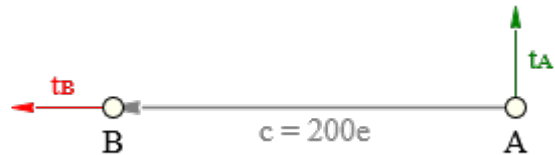

Fig. 14: Single Biarc example.

From these data we extract the vectors: $\quad \mathbf{c} = \begin{pmatrix} -200 \\ 0 \end{pmatrix} e, \quad \mathbf{t}_A = \begin{pmatrix} 0 \\ 1 \end{pmatrix}, \quad \mathbf{t}_B = \begin{pmatrix} -1 \\ 0 \end{pmatrix}. \quad$ Find ...

1. the Biarc angle $\psi$
2. the joint circle radius $R$
3. for the *Equal Chords* Biarc ...
   a. the chord vectors $\mathbf{a}$ and $\mathbf{b}$
   b. the arc angles $\alpha$ and $\beta$
   c. the arc radii $R_A$ and $R_B$
   d. the arc centers $\mathbf{r}_{AA_0}$ and $\mathbf{r}_{BB_0}$
4. for the *Cubic Midpoint* Biarc ...
   a. the chord vectors $\mathbf{a}$ and $\mathbf{b}$
   b. the control points $\mathbf{r}_{AA_1}$ and $\mathbf{r}_{BB_1}$
   c. the arc angles $\alpha$ and $\beta$
   d. the arc radii $R_A$ and $R_B$

**1. The Biarc angle $\psi$**

The characteristic Biarc angle $\psi$ can be seen directly in Figure 14 or determined using equation (1). It is the directed angle from $\mathbf{t}_A$ to $\mathbf{t}_B$.

$$\psi = \tan^{-1} \frac{\tilde{\mathbf{t}}_A \mathbf{t}_B}{\mathbf{t}_A \mathbf{t}_B} = \tan^{-1} \frac{\begin{pmatrix} -1 \\ 0 \end{pmatrix} \begin{pmatrix} -1 \\ 0 \end{pmatrix}}{\begin{pmatrix} 0 \\ 1 \end{pmatrix} \begin{pmatrix} -1 \\ 0 \end{pmatrix}} = 90°$$

**2. The joint circle radius $R$**

is found via equation (3)

$$R = \frac{c}{2 \sin \frac{\psi}{2}} = \frac{200e}{2 \sin 45°} = \underline{141.4e}$$

**3) The *Equal Chords* Biarc**

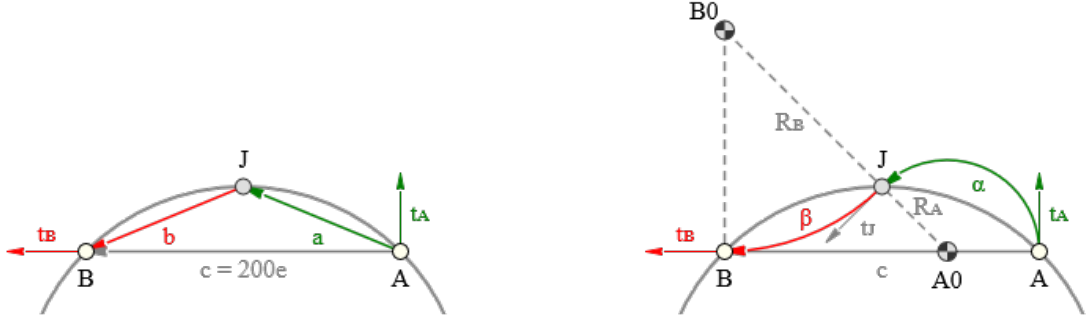

Fig. 15: Equal Chords Biarc.

*3.a) The chord vectors **a** and **b***

According to section 8.1 we get the chord vector $\mathbf{a} = \mathbf{r}_{AJ}$ for the equal chords Biarc by setting $u = 0$ in equation (10) or directly take equation (8.1).

$$\mathbf{a} = \mathbf{r}_{AJ_0} = \tfrac{1}{2}\left(\mathbf{c} - \tan\tfrac{\psi}{4}\tilde{\mathbf{c}}\right) = \tfrac{1}{2}\left[\begin{pmatrix}-200\\0\end{pmatrix}e - \tan 22.5° \begin{pmatrix}0\\-200\end{pmatrix}e\right] = \underline{\begin{pmatrix}-100\\41.4\end{pmatrix}e}$$

Chord vector **b** is easy to get via equation (5)

$$\mathbf{b} = \mathbf{c} - \mathbf{a} = \begin{pmatrix}-200\\0\end{pmatrix}e - \begin{pmatrix}-100\\41.4\end{pmatrix}e = \underline{\begin{pmatrix}-100\\-41.4\end{pmatrix}e}$$

*3.b) The arc angles $\alpha$ and $\beta$*

Arc angles $\alpha$ and $\beta$ are calculated by equations (6) and (2)

$$\alpha = 2\tan^{-1}\frac{\tilde{\mathbf{t}}_A \mathbf{a}}{\mathbf{t}_A \mathbf{a}} = 2\tan^{-1}\frac{\begin{pmatrix}-1\\0\end{pmatrix}\begin{pmatrix}-100\\41.4\end{pmatrix}}{\begin{pmatrix}0\\1\end{pmatrix}\begin{pmatrix}-100\\41.4\end{pmatrix}} = 2\tan^{-1}\frac{100}{41.4} = \underline{135°}$$

$$\beta = \psi - \alpha = 90° - 135° = \underline{-45°}$$

*3.c) The arc radii $R_A$ and $R_B$*

Arc radii are calculated by equation (7). Note that we have signed radii.

$$R_A = \frac{a}{2\sin\frac{\alpha}{2}} = \frac{108.2e}{2\sin 67.5°} = \underline{58.6e} \quad \text{and} \quad R_B = \frac{b}{2\sin\frac{\beta}{2}} = \frac{108.2e}{2\sin -22.5°} = \underline{-141.4e}$$

*3.d) The arc centers $\mathbf{r}_{AA_0}$ and $\mathbf{r}_{BB_0}$*

We use equation (7.1) for calculating the arc center locations.

$$\mathbf{r}_{AA_0} = R_A \tilde{\mathbf{t}}_A = 58.6e \begin{pmatrix}-1\\0\end{pmatrix} = \underline{\begin{pmatrix}-58.6\\0\end{pmatrix}e}$$

$$\mathbf{r}_{BB_0} = R_B \tilde{\mathbf{t}}_B = -141.4e \begin{pmatrix}0\\-1\end{pmatrix} = \underline{\begin{pmatrix}0\\141.4\end{pmatrix}e}$$

## 4) The *Cubic Midpoint* Biarc

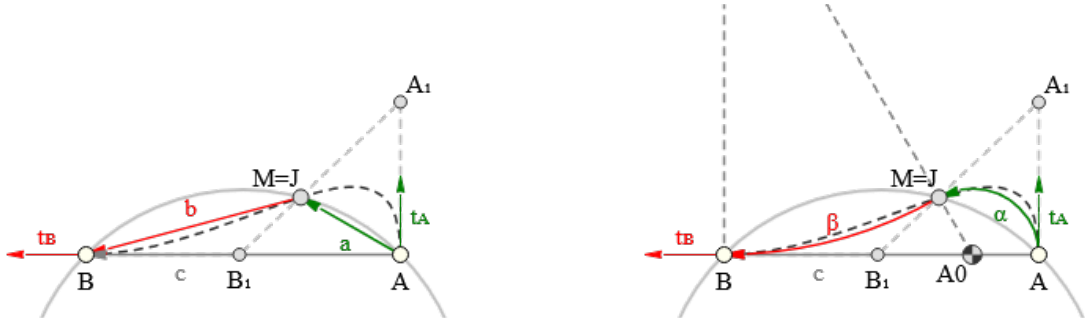

Fig. 16: Cubic Midpoint Biarc.

*4.a) The chord vectors* **a** *and* **b**

To determine the midpoint of the cubic Bezier curve M we need to get some helper variables first (see Theorem 10).

$$\mathbf{t}_{BA} = \mathbf{t}_A - \mathbf{t}_B = \begin{pmatrix} 1 \\ 1 \end{pmatrix}, \quad \kappa = \frac{\tilde{\mathbf{t}}_{BA}\mathbf{c}}{\mathbf{t}_{BA}^2} = 100e, \quad \tan\tfrac{\psi}{2} = 1, \quad \sin^2\tfrac{\psi}{2} = \tfrac{1}{2}$$

$$h = \frac{-4\kappa}{3\tan\tfrac{\psi}{2}} + \sqrt{\frac{16\kappa^2}{9\tan^2\tfrac{\psi}{2}} + \frac{4c^2}{9\sin^2\tfrac{\psi}{2}}} = \ldots = 97.6e$$

Now from equation (18) we get chord vector **a**.

$$\mathbf{a} = \tfrac{1}{2}\mathbf{c} + \tfrac{3}{8}h\mathbf{t}_{BA} = \tfrac{1}{2}\begin{pmatrix} -200 \\ 0 \end{pmatrix} + \tfrac{3}{8}97.6e\begin{pmatrix} 1 \\ 1 \end{pmatrix} = \underline{\begin{pmatrix} -63.4 \\ 36.6 \end{pmatrix} e}$$

and chord vector **b**.

$$\mathbf{b} = \mathbf{c} - \mathbf{a} = \begin{pmatrix} -200 \\ 0 \end{pmatrix} e - \begin{pmatrix} -63.4 \\ 36.6 \end{pmatrix} e = \underline{\begin{pmatrix} -136.6 \\ -36.6 \end{pmatrix} e}$$

*4.b) The control points* $\mathbf{A}_1$ *and* $\mathbf{B}_1$ *(not needed for further calculation, just for illustration)*

The relative control point locations are (see Theorem 10)

$$\mathbf{r}_{AA_1} = h\mathbf{t}_A = 97.6e\begin{pmatrix} 0 \\ 1 \end{pmatrix} = \underline{\begin{pmatrix} 0 \\ 97.6 \end{pmatrix} e}, \quad \mathbf{r}_{BB_1} = -h\mathbf{t}_B = -97.6e\begin{pmatrix} -1 \\ 0 \end{pmatrix} = \underline{\begin{pmatrix} 97.6 \\ 0 \end{pmatrix} e}.$$

*4.c) The arc angles* $\alpha$ *and* $\beta$

According to 3.c) arc angles $\alpha$ and $\beta$ are calculated by equations (6) and (2)

$$\alpha = 2\tan^{-1}\frac{\tilde{\mathbf{t}}_A \mathbf{a}}{\mathbf{t}_A \mathbf{a}} = 2\tan^{-1}\frac{\begin{pmatrix} -1 \\ 0 \end{pmatrix}\begin{pmatrix} -63.4 \\ 36.6 \end{pmatrix}}{\begin{pmatrix} 0 \\ 1 \end{pmatrix}\begin{pmatrix} -63.4 \\ 36.6 \end{pmatrix}} = 2\tan^{-1}\frac{63.4}{36.6} = \underline{120°}$$

$$\beta = \psi - \alpha = 90° - 120° = \underline{-30°}$$

*4.d) The arc radii* $R_A$ *and* $R_B$

Arc radii are calculated by equation (7) again.

$$R_A = \frac{a}{2\sin\tfrac{\alpha}{2}} = \frac{73.2e}{2\sin 60°} = \underline{42.3e} \quad \text{and} \quad R_B = \frac{b}{2\sin\tfrac{\beta}{2}} = \frac{141.4e}{2\sin -15°} = \underline{-273.2e}$$

# 7. Conclusion

Biarcs are often used to approximate curves via arc splines having $G^1$ continuity, as required in CNC machining and robot path planning, for example. For a pair of oriented $G^1$ Hermite data, two arc segments – named *biarc* – are generally required to hit both end points and their associated tangents.

This article takes a purely geometric view of biarcs from the perspective of symplectic geometry. For two Hermite data sets, the corresponding biarc is not uniquely defined, but the join point of both arcs belongs to a one-parameter family of points, all of which lie on a circle.

The fundamental properties of a biarc and its joint circle are presented. Based on these insights, various well-known strategies from the biarc literature for selecting an advantageous distinct join point from its one-parameter family are discussed. In addition a new proposal is presented that uses the center point of a cubic spline function.

Some illustrative examples show the influence of the strategy used to select the join point on the resulting arc spline curve.

The use of symplectic geometry in $\mathbb{R}^2$ results in a set of new pure vector equations, which can be easily reused in engineering software and computer graphics.

# References

[1] Arnold, V.I., "Mathematical methods in Classical Mechanics", Springer (1978).

[2] Arnold, V.I., Givental, A.B., "Symplectic Geometry",
https://www.maths.ed.ac.uk/~v1ranick/papers/arnogive.pdf

[3] Gotay, M.J., Isenberg, J.A., "The symplectification of science",
https://www.pims.math.ca/~gotay/Symplectization(E).pdf

[4] Gössner, S. "Symplectic Geometry for Engineers - Fundamentals", Dec. 2019;
https://arxiv.org/abs/2404.01091

[5] Goessner, S., "Symplectification of Circular Arcs and Arc Splines" Aug. 2025,
https://arxiv.org/abs/2508.07726

[6] Bolton, K.M., "Biarc curves", Computer Aided Design 7 (1975).

[7] Sabin, M.A., "The use of piecewise forms for the numerical representation of shape",
Report 60/1977, Computer and Automation Institute, Hungary.

[8] Chandrupatla, T.R., Osler, T.J., "Planar Biarc Curves – A Geometric View", (2011)
https://www.researchgate.net/publication/265802980_Planar_Biarc_Curves_-_A_Geometric_View

[9] Meek, D.S., Walton, D.J. "Approximating smooth planar curves by arc splines", J. Comput. Appl. Math. 59 (1995).

[10] Meek, D.S., Walton, D.J. "The family of biarcs that matches planar, two-point G1 Hermite data", J. Comput. Appl. Math. 212 (2008).

[11] Šír, Z., Feichtinger, R., Jüttler, B., "Approximating curves and their offsets using biarcs and Pythagorean hodograph quintics", Computer-Aided Design, Volume 38, 2006

[12] Yong, J.H., Hu, S.M., Sun, J.G., "A note on approximation of discrete databy G1 arc splines", Computer-Aided Design, Volume 31, 1999

[13] Hoschek, J., Lasser, D., "Grundlagen der geometrischen Datenverarbeitung", Teubner (1992).

[14] "Centripetal Catmull–Rom spline", https://en.wikipedia.org/wiki/Centripetal_Catmull–Rom_spline